\documentclass{article}
\usepackage[ansinew]{inputenc}
\usepackage[dvips]{graphicx}
\usepackage[T1]{fontenc}
\usepackage{amsmath,enumerate,amssymb,hhline,amsthm}

 \newtheorem{thm}{Theorem}[section]
 \newtheorem{cor}[thm]{Corollary}
 \newtheorem{lem}[thm]{Lemma}
 
 \newtheorem{conj}[thm]{Conjecture} 
 \theoremstyle{definition}
 \newtheorem{rem}[thm]{Remark}

\begin{document}

\title{On the shifted convolution problem in mean}

\author{Eeva Suvitie}

\date{}

\maketitle

\let\thefootnote\relax\footnotetext{This
research was supported by the Finnish Cultural Foundation and by the project 138522 of the Academy of Finland, and
accomplished while the author was visiting \'{E}cole Polytechnique F\'{e}d\'{e}rale de Lausanne for
a time period of one year between April 2011 - April 2012.}

\begin{abstract}
We study the following mean value of the shifted convolution problem:
\[
\sum_{f\sim F}\sum_{n\sim N}\left|\sum_{l\sim L}t(n+l)t(n+l+f)
\right|^{2},
\]
over the Hecke eigenvalues of a fixed non-holomorphic cusp form with quantities $N\geq 1$, $1\leq L\leq N^{1-\varepsilon}$ and $1\leq F\ll N^{2/5}$. We attain a result also for a weighted case. 
Furthermore, we point out that the proof yields 
analogous upper bounds for the shifted convolution problem over the Fourier coefficients of a fixed holomorphic cusp form in mean. \\ \\
\textbf{Mathematics Subject Classification (2010):} Primary 11F30; Secondary 11F72, 11P55. \\ \textbf{Keywords:} The shifted convolution problem, circle method, spectral theory.
\end{abstract}

\section{Introduction}

The additive divisor problem and its analogs for the Fourier coefficients of cusp forms, together called the shifted
convolution problem, have been under thorough study since the early 20th century.
A brief survey of the history is given in our previous paper \cite{Suvitie2}, leaning heavily on works by Blomer and Harcos \cite{Harcos} and Motohashi \cite{Motohashi4}. 

In \cite{Suvitie2} we studied this problem in a mean value sense, the motivation arising from Jutila's paper \cite{Jutila9}, where he examined the shifted convolution problem
over the Fourier coefficients of a holomorphic cusp form of weight $k$ for the full modular group. He required an estimate for the mean value
\[
\sum_{0\leq f\leq F}\bigg|\sum_{1\leq n\leq N}
a(n)\overline{a(n+f)}\bigg|^{2}
\]
over the same Fourier coefficients (Lemma 3), a similar argument yielding also the following estimate for a triple sum:
\begin{equation}\label{aasumma}
\sum_{0\leq f\leq F}\sum_{1\leq n\leq N}\bigg|\sum_{1\leq l\leq L}
a(n+l)\overline{a(n+f+l)}\bigg|^{2}\ll (N+F)^{k}N^{k}L
\end{equation}
for all $N,F\geq 1$ and $1\leq L\leq N$.
Instantaneously the proof gives also an upper bound 
for the analogous sum over the Hecke eigenvalues of a non-holomorphic cusp form in mean, as stated in Lemma 6 of \cite{Jutila7} :  For $N,F\geq 1$ and $1\leq L\leq N$ 
\begin{equation}\label{Tsumma}
\sum_{0\leq f\leq F}\sum_{1\leq n\leq N}\bigg|\sum_{1\leq l\leq L}
t(n+l)t(n+f+l)\bigg|^ {2}\ll (N+F)^{1+\varepsilon}NL.
\end{equation}
The novelty of the estimation of the triple sum lies
in its sensitiveness for the length of the innermost sum over $l$.
A similar result, sensitive for the size of the shift, is needed in the doctoral thesis of the author (\cite{Suvitie}, Lemma 3.4). 
Also, an analogous result in a more general setting appears on p. 81 in the paper \cite{Blomer} by Blomer, Harcos and Michel, leading
to the same bounds \eqref{aasumma}, \eqref{Tsumma}.
In the three papers \cite{Blomer}, \cite{Jutila9}, \cite{Jutila7} the motivation underlying the study of the mean values, including either double or triple sums, has been to obtain information about the 
upper bounds of the original shifted convolution sums, while on the other hand the author needed her estimate in \cite{Suvitie} to
estimate a certain spectral sum over inner products involving a holomorphic cusp form and Maass forms.

However, applying the methods of the earlier proofs seemed problematic when it came to extending the results to the analogous case of the additive
divisor problem.
In our paper \cite{Suvitie2} we were able to prove the following result:
Let $N\geq 1$, $1\leq L\leq N$ and $1\leq F\ll N^{1-\varepsilon}$. Then
\[
\sum_{f\sim F}\sum_{n\sim N}\left|\sum_{l\sim L}d(n+l)d(n+l+f)-\frac{6}{\pi^{2}}\int_{(L+n)/f}^{(2L+n)/f}m(x;f)dx
\right|^{2}
\]
\begin{equation}\label{paatermi}
\ll N^{2+\varepsilon}+N^{1+\varepsilon}LF.
\end{equation}
Here $m(x;f)$ is as in (1.12) in \cite{Motohashi4}, and
the notation $m\sim M$ stands for $M< m\leq2M$. The proof was based on a spectral decomposition of the
shifted convolution sum over the divisor function (see \cite{Motohashi4}, Theorem 3).
Following the analogous argument for the case of a holomorphic cusp form we also immediately attained
\begin{equation}\label{edellinena}
\sum_{f\sim F}\sum_{n\sim N}\left|\sum_{l\sim L}a(n+l)\overline{a(n+l+f)}
\right|^{2}
\ll
N^{2k+\varepsilon}+N^{2k-1+\varepsilon}LF
\end{equation}
for $N\geq 1$, $1\leq L\leq N$ and $1\leq F\ll N^{1-\varepsilon}$. 
In the case of a non-holomorphic cusp form we faced the problem of lacking a proper analogy for the spectral decomposition of the shifted convolution sum
in question, crucial for our proof. In this case a spectral decomposition exists, but there appears also an arithmetic correction term, which leads to additional problems. 
Therefore we just stated the following conjecture:
\begin{conj}\label{con}
Let $N\geq 1$, $1\leq L\leq N$ and $1\leq F\ll N^{1-\varepsilon}$. Then
\[
\sum_{f\sim F}\sum_{n\sim N}\left|\sum_{l\sim L}t(n+l)t(n+l+f)
\right|^{2}
\ll
N^{2+\varepsilon}+N^{1+\varepsilon}LF.
\]
\end{conj}
In this paper we attack this conjecture through Jutila's version of the classical circle method, aiming at a unified bound for all three cases of
triple sums.

\section{Results}

The following notation will be adopted:  
Vinogradov's relation $f(z)\ll g(z)$ is another notation for
$f(z)=\mathcal{O}(g(z))$. Let $\varepsilon$ stand generally for a small positive number,
not necessarily the same at each occurrence. 

Jutila's version of the circle method was originally introduced in \cite{Jutila9}, based on a well-distribution property for
rational numbers examined in \cite{Jutila11}. It was further generalized in
\cite{Jutila12}-\cite{Jutila10}, the last two standing out as our primary guides in the course of our proof.
As in \cite{Jutila6} we interpret our sum inside the square as a Fourier coefficient of a suitable continuous function $\psi$ of
period $1$;
\begin{equation}\label{aloitus}
b_{f}= \int_{0}^{1}\psi(x)e(-f x)\,dx, \;\; e(\alpha)=\exp(2\pi i \alpha),
\end{equation}
that we estimate through an approximation constructed as follows. Let $Q>0$ be a large parameter and $w(q)$ stand for a function such that $0\leq w(q)\leq 1$ for all natural numbers $q$.
Moreover let $w(q)=0$ for
$q\notin [Q,2Q]$ and choose $w$ so that
\[
\Lambda=\sum_{q=1}^{\infty}w(q)\phi(q)\gg Q^{2}
\]
with $\phi$ standing for Euler's totient function.
Let $\Delta\in (0,1/3)$ be another parameter, and let $0\leq \nu(x)\ll 1$ stand for a piecewise monotonic continuous function supported in the interval $[-\Delta,-\Delta/2]$,
satisfying
\[
\frac{1}{2\Delta}\int_{-\infty}^{\infty}\nu(x)\,dx=1.
\]
Note that the asymmetric choice of the support will be convenient later.
We denote $2\Delta \Lambda$ by $\lambda$. Then we define an approximation for the characteristic function $\chi(x)$ of the interval $[0,1]$;
\[
\chi^{\ast}(x)=\frac{1}{\lambda}\sum_{q=1}^{\infty}w(q)\sum_{\substack{a=1 \\ (a,q)=1}}^{q}\nu^{\ast}\left(\frac{a}{q}-x \right),
\]
where $\nu^{\ast}(x)$ is the extension of the function $\nu(x)$ to a periodic function with period $1$. Now
\begin{equation}\label{khii}
\chi^{\ast}(x)=1+\frac{1}{\lambda}\sum_{d=1}^{\infty}d\sum_{m\neq 0}a_{dm}e(-dmx)\sum_{r=1}^{\infty}w(dr)\mu(r),
\end{equation}
where the Fourier coefficients of $\nu^{\ast}$ are denoted by $a_{\beta}$ and $\mu(r)$ stands for the M\"obius function.
If the function $\nu(x)$ is chosen sufficiently smooth, we have
\begin{equation}\label{aFourier}
a_{\beta}\ll (\beta^{A+1}\Delta^{A}+\Delta^{-1})^{-1}
\end{equation} 
for any large constant $A>0$.
Now an approximation to $b_{f}$ is
\begin{equation}\label{beexii}
b_{f}^{\ast}= \int_{0}^{1}\chi^{\ast}(x)\psi(x)e(-f x)\,dx,
\end{equation}
whence we get
\begin{equation}\label{tahti}
b_{f}^{\ast}=\frac{1}{\lambda}\sum_{q=1}^{\infty}w(q)\sum_{\substack{a=1 \\ (a,q)=1}}^{q}\int_{-\infty}^{\infty}\nu\left(\frac{a}{q}-x \right)\psi(x)e(-f x)\, dx.
\end{equation}
From Theorem 1 in \cite{Jutila12} we have the estimate
\begin{equation}\label{vee}
V=\int_{0}^{1}\left(\chi^{\ast}(x)-\chi(x)\right)^{2}\, dx\ll \lambda^{-1}\log^{3}(\Delta^{-1});
\end{equation}
note that our conditions on $\nu(x)$ above suffice, although differing from those in \cite{Jutila12}.

In \cite{Jutila6} Jutila utilized the bound (Theorem 1, \cite{Jutila6})
\[
\sum_{f=0}^{F}\beta_{f}b_{f}\ll \|\psi\|_{\infty}Q^{-A}+\max_{|\xi|\leq Q^{1+\varepsilon}}\left|\sum_{f=0}^{F}\beta_{f}b_{\xi+f}^{\ast}\right|
\]
for the study of the shifted convolution problem over the Fourier coefficients of a non-holomorphic cusp form. 
Here  $A,\varepsilon>0$ are any fixed constants, $0\leq F\ll Q^{c}$ for some constant $c>0$, $\Delta\asymp Q^{-1}$ and the $\beta_{f}$ are arbitrary 
complex numbers such that $|\beta_{f}|\leq 1$.
In our case, however, the presence of the double sum and the square produces new problems and we need again a refined estimate for the difference
$b_{f}-b_{f}^{\ast}$, presented in Lemma \ref{uusiarvio} below.

We follow an unpublished preprint by Jutila \cite{Jutila10}, where the following deduction is made:
We write \eqref{khii} as
\[
\chi^{\ast}(x)=1-E(x),
\]
where
\[
E(x)=\sum_{\xi\neq 0}c_{\xi}e(-\xi x)
\]
with
\begin{equation}\label{seean}
c_{\xi}=-\frac{a_{\xi}}{\lambda}\sum_{d|\xi}d\sum_{r=1}^{\infty}w(dr)\mu(r)\ll \lambda^{-1}Q\, |a_{\xi}|\, d(\xi).
\end{equation}
We note that for any non-negative integer $k$ and $0\leq x\leq 1$
\begin{equation}\label{Esumma}
\chi(x)=1=\chi^{\ast}(x)(1+E(x)+\ldots+E^{k}(x))+E^{k+1}(x).
\end{equation}
Now if we choose $\Delta=Q^{-1+\delta_{2}}$ for a small constant $\delta_{2}>0$ (unlike in \cite{Jutila6}), then by \eqref{aFourier}, \eqref{seean} and
straightforward estimates
\begin{equation}\label{Earvio}
E(x)=1-\chi^{\ast}(x)\ll Q^{-\varepsilon}.
\end{equation}
Therefore, by \eqref{aloitus} and \eqref{Esumma}, when we choose $k$ large enough,
\[
b_{f}=b_{f}^{\ast}+\sum_{\xi=-\infty}^{\infty}d_{\xi}b_{f+\xi}^{\ast}+\mathcal{O}(\|\psi\|_{\infty}Q^{-A}),
\]
where
\[
d_{\xi}=\sum_{\substack{\xi_{1},\xi_{2}\neq 0 \\ \xi_{1}+\xi_{2}=\xi}}c_{\xi_{1}}c_{\xi_{2}}+\cdots+\sum_{\substack{\xi_{1},\ldots ,\xi_{k}\neq 0 \\ \xi_{1}+\cdots +\xi_{k}=\xi}}c_{\xi_{1}}\cdots c_{\xi_{k}}
+\bigg\{\begin{array}{ll} 0, & \xi=0,\\
c_{\xi}, & \xi\neq 0.\\
\end{array}
\]
On the other hand, the $\xi$th Fourier coefficient of $E^{k}(x)$, $k\geq 2$, is by \eqref{Earvio}
\[
\int_{0}^{1}E^{k}(x)e(-\xi x)\, dx\ll\int_{0}^{1}|E^{2}(x)|\, dx=V,
\]
and hence by \eqref{vee} and \eqref{seean} always $d_{\xi}\ll Q^{-1}d(\xi)$. For large values of $\xi$ even more can be achieved, and we may truncate the $\xi$-sum above
by using \eqref{aFourier}, \eqref{beexii} and the upper bound \eqref{seean} for $c_{\xi}$, attaining
\begin{lem}[Jutila]\label{uusiarvio}
For all $f\in \mathbb{Z}$ 
\[
b_{f}=b_{f}^{\ast}+\sum_{|\xi|\ll Q}d_{\xi}b_{f+\xi}^{\ast}+\mathcal{O}(\| \psi\|_{\infty}Q^{-A}),
\]
with $d_{\xi}\ll Q^{-1}d(\xi)$.
\end{lem}

\vspace{0,5cm}

In our study we shall first consider a weighted sum: We always let $N\geq 1$, $L\leq N^{1-\varepsilon}$. Furthermore let $L^{-1+\varepsilon}\leq \delta\leq 1/4$ and 
$0\leq W_{n}(x)\leq 1$ be a real-valued smooth weight function supported on some interval 
$[BL+n,CL+n]$, where $B,C$ are independent of the variable $n$, $1\leq B\leq C\leq 2$ and $C-B\ll \delta$. We suppose
$W_{n}^{(\nu)}(x)\ll_{\nu}(\delta L)^{-\nu}$ for each $\nu\geq 0$ and $x\in \mathbb{R}$. Moreover let $0\leq W_{n}^{0}(x)\leq 1$ be a real-valued smooth weight function supported on some interval 
$[B'L+n,C'L+n]$, where also $B',C'$ are independent of the variable $n$, $1/2 \leq B'\leq C'\ll 1$, $C'-B'\ll \delta$, and
$(W_{n}^{0})^{(\nu)}(x)\ll_{\nu}(\delta L)^{-\nu}$ for each $\nu\geq 0$ and $x\in \mathbb{R}$. We choose $Q=\delta LN^{-\delta_{1}}$ for some fixed small $\delta_{1}>0$, and assume
at first that $L\gg N^{\varepsilon}$ so that always $Q\gg N^{\varepsilon}$. Moreover, let $\Delta=Q^{-1+\delta_{1}/2}$, that is, we specify $\delta_{2}=\delta_{1}/2$ above. 
Finally, let $\psi_{n}(x)$ stand for 
the periodic function $S(W_{n}^{0},x)S(W_{n},-x)$ with
\[
S(W,x)=\sum_{m=1}^{\infty}W(m)t(m)e(mx).
\]
Now $\psi_{n}(x)$ depends on the variable $n$, so we indicate this dependence also in its $f$th Fourier coefficient 
\[
b_{f}(n)=\sum_{m=1}^{\infty}t(m)t(m+f)W_{n}(m)W_{n}^{0}(m+f).
\]

Now we proceed to present our results; the proofs are given in Section 4.
First, for the estimation of the weighted sum we have
\begin{lem}\label{tarkein}
Let $1\leq F\ll \delta L$. Then
\[
\sum_{|f|\sim F}\sum_{n\sim N}\left |b_{f}^{\ast}(n)
\right|^{2} \ll N^{3+\varepsilon}(\delta L)^{-2}+N^{2+\varepsilon}+N^{1+\varepsilon}\delta LF.
\]
\end{lem}

In case $f=0$ we get
\begin{lem}\label{nolla} 
\[
\sum_{n\sim N}\left |b_{0}^{\ast}(n)
\right|^{2} \ll N^{1+\varepsilon}(\delta L)^{2}.
\]
\end{lem}

Combining the above lemmas, we get the following auxiliary result:
\begin{lem}\label{johtopaatos}
We have
\[
\sum_{n\sim N}\left |b_{f}(n)
\right|^{2} \ll \sum_{n\sim N}\left |b_{f}^{\ast}(n)
\right|^{2}+N^{3+\varepsilon}(\delta L)^{-3}+N^{2+\varepsilon}(\delta L)^{-1}+N^{1+\varepsilon}\delta L
\]
uniformly for $1\leq f\ll \delta L$.
\end{lem}

Finally we have an upper bound for the weighted sum:
\begin{thm}\label{thm1} 
Let $N\geq 1$, $1\leq L\leq N^{1-\varepsilon}$ and $1\leq F\ll \delta L$, $L^{-1+\varepsilon}\leq \delta\leq 1/4$. Let
$0\leq W_{n}(x)\leq 1$ be a real-valued smooth weight function supported on some interval 
$[BL+n,CL+n]$, where $B$ and $C$ are independent of $n$, $1\leq B\leq C\leq 2$ and $C-B\ll \delta$. Suppose that
$W_{n}^{(\nu)}(x)\ll_{\nu}(\delta L)^{-\nu}$ for each $\nu\geq 0$ and $x\in \mathbb{R}$. Then
\[
\sum_{f\sim F}\sum_{n\sim N}\left|\sum_{l\sim L}t(n+l)t(n+l+f)W_{n}(n+l)
\right|^{2}
\ll
N^{3+\varepsilon}(\delta L)^{-2}+N^{2+\varepsilon}+N^{1+\varepsilon}\delta LF.
\]
\end{thm}

For the non-weighted sum we immediately conclude
\begin{thm}\label{thm2} 
Let $N\geq 1$, $1\leq L\leq N^{1-\varepsilon}$ and $1\leq F\ll N^{2/5}$. Then
\[
\sum_{f\sim F}\sum_{n\sim N}\left|\sum_{l\sim L}t(n+l)t(n+l+f)
\right|^{2}
\ll
N^{2+\varepsilon}F^{1/2}+N^{1+\varepsilon}LF.
\]
\end{thm}
Thus in case $L\leq N^{1-\varepsilon}$, $F\ll \min(N^{2/5},L^{2})$ and $N^{2}\ll L^{4}F$ our result is better than what \eqref{Tsumma} or the trivial estimate $N^{1+\varepsilon}L^{2}F$ yield. 

We reach Conjecture \ref{con} and equivalently the bound \eqref{edellinena} from \cite{Suvitie2} by new methods (see Theorem \ref{thm4} below), with some 
additional restrictions on the parameters:
\begin{cor}\label{cor1}
When $L\ll N^{1-\varepsilon}$, Conjecture \ref{con} holds if $1\ll F\ll N^{2/5}$ and $N^{2}\ll L^{2}F$, or if 
$F\asymp 1$.
\end{cor}
The result is not as strong as in \cite{Suvitie2}, which is a consequence of some abrupt estimates necessary in the course
of our proof.

\vspace{0,5cm}

An analogous argument gives the following result for the case of a holomorphic cusp form:
\begin{thm}\label{thm4} 
Let $N\geq 1$, $1\leq L\leq N^{1-\varepsilon}$ and $1\leq F\ll \delta L$, $L^{-1+\varepsilon}\leq \delta\leq 1/4$. Let the weight function $W_{n}$
be as in Theorem \ref{thm1}. Then
\[
\sum_{f\sim F}\sum_{n\sim N}\left|\sum_{l\sim L}a(n+l)\overline{a(n+l+f)}W_{n}(n+l)
\right|^{2}
\ll
N^{2k+1+\varepsilon}(\delta L)^{-2}+N^{2k+\varepsilon}
\]
\[
+N^{2k-1+\varepsilon}\delta LF.
\]
If $N\geq 1$, $1\leq L\leq N^{1-\varepsilon}$ and $1\leq F\ll N^{2/5}$, then
\[
\sum_{f\sim F}\sum_{n\sim N}\left|\sum_{l\sim L}a(n+l)\overline{a(n+l+f)}
\right|^{2}
\ll
N^{2k+\varepsilon}F^{1/2}+N^{2k-1+\varepsilon}LF.
\]
\end{thm}

\begin{rem}
In case $F\gg \delta L$ we may repeat essentially the same deduction as above, but we run into troubles with the oscillating factor
$e^{ixa(m/p,n)}$ (see \eqref{eeaxaa} below). For example, if $m_{1}<m_{2}$, using the same notation as in the proof of Lemma \ref{tarkein}, and 
having $P\asymp FQ^{\delta_{1}}$, the term $a(m/p,n)$ may be small although its derivative with respect to $p$ is large. It appears that these difficulties may
be overcome by multiplying the exponential factor by some suitable power of $F/\delta L$, and we might still be able to attain an upper bound
better than what \eqref{Tsumma} yields in certain extremal cases, at least when $F\ll N^{1-\varepsilon}$. However, as applying the method in this way is not very elegant or efficient, we omit the case 
$F\gg \delta L$ altogether in this paper, leaving the search for a smoother way of treating it for further study. 
\end{rem}

\begin{rem}
In proving our results above we shall use the earlier estimates \eqref{aasumma} and \eqref{Tsumma}. However, an iterative deduction where our new upper bound would be 
fed back to the proof would not work, as the estimates \eqref{aasumma} and \eqref{Tsumma} are needed especially in the case $F=N$, which is outside
our range for $F$. Also it seems probable that in this case the old bounds are the best that we can hope to achieve.
\end{rem}

\begin{rem}
This time our proof cannot be easily extended to the case of the divisor function $d(n)$, because of the problematic main term appearing in 
\eqref{paatermi}, so the quest for a unified method for attacking all three analogous cases still remains.
\end{rem}

\section{Needed notation and auxiliary lemmas}

\subsection{Cusp forms}

We recall the notation and some results on the cusp forms appearing in our paper.
For proofs and general reference the reader is referred to Motohashi's monograph \cite{Motohashi2}.

We confine ourselves to cusp forms for the full modular group
$\Gamma=SL_{2}(\mathbb{Z})$ operating through M\"{o}bius
transformations on the upper half-plane
$\mathbb{H}$. A
\textit{holomorphic cusp form}
$F(z):\mathbb{H}\rightarrow\mathbb{C}$ of weight $k\in\mathbb{Z}$
with respect to $\Gamma$ can be represented by its Fourier series
\[
F(z)=\sum^{\infty}_{n=1}a(n)e(nz).
\]
We may assume that $k$ is even and $k\geq 12$, otherwise $F(z)$ is
trivial.
Let
\[
\{\psi_{j,k}\;|\; 1\leq j\leq\vartheta(k)\}
\]
be an orthonormal basis of the unitary space of holomorphic cusp
forms of weight $k$, and write
\[
\psi_{j,k}(z)=\sum_{n=1}^{\infty}\rho_{j,k}(n)n^{\frac{k-1}{2}}e(nz).
\]
We may suppose that the basis vectors are eigenfunctions of the
Hecke operators $T_{k}(n)$ for all positive integers $n$. Thus, in particular,
$T_{k}(n)\psi_{j,k}=t_{j,k}(n)\psi_{j,k}$ for certain real numbers
$t_{j,k}(n)$, which we call Hecke eigenvalues. Comparing the Fourier
coefficients on both sides, one may verify that
$\rho_{j,k}(n)=\rho_{j,k}(1)t_{j,k}(n)$ for all $n\geq 1$, $1\leq
j\leq \vartheta(k)$.
We put 
\[
a_{k}=2^{2-2k}\pi^{-k-1}(k-1)!,
\]
whence by Eq. (2.2.10) in \cite{Motohashi2}
\begin{equation}\label{aakoo}
a_{k}\sum_{j=1}^{\vartheta(k)}|\rho_{j,k}(m)|^{2}\ll km^{1/2+\varepsilon}
\end{equation}
for any integers $k,m\geq 1$.
Furthermore we recall the bound
\begin{equation}\label{tee}
|t_{j,k}(n)|\leq d(n)\ll n^{\varepsilon}
\end{equation}
by Deligne \cite{Deligne}.

A \textit{non-holomorphic cusp form} $u(z)=u(x+iy):\mathbb{H}\rightarrow\mathbb{C}$ is a
non-constant real-analytic $\Gamma$-invariant function in the
upper half-plane, square-integrable with respect to the hyperbolic
measure $d\mu(z)=\frac{dx\,dy}{y^{2}}$ over a fundamental domain
of $\Gamma$. Also $u(z)$ is an eigenfunction of the non-euclidean Laplacian
 $\Delta=-y^2(\frac{\partial^2}{\partial x^2}+\frac{\partial^2}{\partial
 y^2})$, and the corresponding eigenvalue can be written
as $1/4+\kappa^2$ with $\kappa>0$.
The Fourier series expansion for $u(z)$ is then of the form
\[
u(z)=y^{1/2}\sum_{n\neq 0}\rho(n)K_{i\kappa}(2\pi |n|y)e(nx)
\]
with $K_{\nu}$ a Bessel function of imaginary argument. We may
suppose that our cusp forms are eigenfunctions of the Hecke
operators $T(n)$ for all positive integers $n$ and that $u(x+iy)$
is even or odd as a function of $x$. Thus $T(n)u=t(n)u$ for
certain real numbers $t(n)$, which are again called Hecke eigenvalues,
and $u(-\overline{z})=\pm u(z)$. Comparing the Fourier coefficients on
both sides, one may verify that for all $n\geq 1$, $\rho(n)=\rho(1)t(n)$ and
$\rho(-n)=\varepsilon \rho(n)$, with $\varepsilon=\pm 1$ the parity sign of the cusp form in question. 
The \textit{Maass (wave) forms} $u_{j}$
constitute an orthonormal set of non-holomorphic cusp forms
arranged so that the corresponding parameters $\kappa_{j}$
determined by the eigenvalues $1/4+\kappa_{j}^{2}$ lie in an
increasing order. We write $\rho_{j}(n)$ and $t_{j}(n)$ for the corresponding Fourier coefficients and Hecke
eigenvalues.
We let
\[
\alpha_{j}=|\rho_{j}(1)|^{2}/\cosh(\pi\kappa_{j}).
\]

As counterparts for \eqref{tee} we have 
the following two estimates for all $\varepsilon>0$ and $N\geq 1$: First the classical one due to
Iwaniec \cite{Iwaniec4}, Lemma 1;
\[
\sum_{n\leq N}t^{2}(n)\ll \kappa^{\varepsilon}N.
\]
Further it is known that
\begin{equation}\label{Hoff}
\sum_{n\leq N}t^{4}(n)\ll \kappa^{\varepsilon} N^{1+\varepsilon}
\end{equation}
by \cite{Hoffstein}, Lemma 2.1.

\subsection{Some useful lemmas}

In this section we shall gather some auxiliary results.

For the Riemann zeta function $\zeta$ it is known that
\[
\frac{1}{\zeta(1+it)}\ll \log |t|,
\]
as $|t|\geq 1$. For a proof, see e.g. 
\cite{Titchmarsh}, p. 132. 

We
have an important tool arising from spectral theory:
\begin{lem}[The spectral large sieve]
For $K\geq 1$, $1\leq \Delta\leq K$, $M\geq 1$ and any complex
numbers $a_{m}$ we have
\[
\sum_{K\leq\kappa_{j}\leq K+\Delta}\alpha_{j}\left|\sum_{m\leq
M}a_{m}t_{j}(m)\right|^{2}\ll
(K\Delta+M)(KM)^{\varepsilon}\sum_{m\leq M}|a_{m}|^{2}.
\]
\end{lem}
For a proof, see Theorem 1.1 in \cite{Jutila4} or Theorem 3.3 in
\cite{Motohashi2}.
The continuous analogue is the following:
\begin{lem}\label{jatkuvasss}
For $K$ real, $\Delta\geq 0$, $M\geq 1$ and any complex
numbers $a_{m}$ we have
\[
\int_{K}^{K+\Delta}\left|\sum_{m\leq
M}a_{m}\sigma_{2ir}(m)m^{-ir}\right|^{2}\, dr\ll
(\Delta^{2}+M)M^{\varepsilon}\sum_{m\leq M}|a_{m}|^{2}
\]
uniformly in $K$, as $\sigma_{2ir}(m)=\sum_{d|m}d^{2ir}$.
\end{lem}
A proof can be found in \cite{Suvitie}, Lemma 1.12.

Next we introduce a Vorono\u\i \, type summation formula:
\begin{lem}\label{Voronoi}
Let $W$ be a smooth function with compact support in $(0,\infty)$, and $a$, $q$ integers such that $q>0$ and $(a,q)=1$.  Then
\[
\sum_{m=1}^{\infty}t(m)e\left(\frac{ma}{q}\right)W(m)=\frac{\pi i}{q\sinh(\pi \kappa)}
\sum_{m=1}^{\infty}t(m)e\left(-\frac{m\overline{a}}{q}\right)
\]
\[
\times
\int_{0}^{\infty}\left(J_{2i\kappa}\left(\frac{4\pi\sqrt{mx}}{q}\right)-J_{-2i\kappa}\left(\frac{4\pi\sqrt{mx}}{q}\right)\right)
W(x)\, dx+\frac{4\varepsilon\cosh(\pi\kappa)}{q}
\]
\[
\times\sum_{m=1}^{\infty}t(m)e\left(\frac{m\overline{a}}{q}\right)\int_{0}^{\infty}K_{2i\kappa}\left(\frac{4\pi\sqrt{mx}}{q}\right)
W(x)\, dx,
\]
where $a\overline{a}\equiv 1$ (mod $q$), $\varepsilon$ is the parity sign of the non-holomorphic cusp form in question 
and $\kappa$ is the related spectral parameter. $J_{\nu}$ and $K_{\nu}$ are the Bessel functions of the first kind and of imaginary argument, respectively.
\end{lem}
For a proof, see \cite{Meurman2}, Theorem 2.

Another crucial lemma is the following identity due to Kuznetsov involving the Kloosterman sums 
\[
S(m,n;q)=\sum_{\substack{a=1\\ (a,q)=1}}^{q}e\left(\frac{ma+n\overline{a}}{q}\right):
\]
\begin{lem}\label{Kuzn}
Let $\psi(x)$ be a smooth function with compact support in $(0,\infty)$ and $m,n$ be positive integers. Then
\[
\sum_{q=1}^{\infty}q^{-1}S(m,n;q)\psi\left(\frac{4\pi\sqrt{mn}}{q}\right)=\sum_{j=1}^{\infty}\alpha_{j}t_{j}(m)t_{j}(n)\hat{\psi}(\kappa_{j})
\]
\begin{equation}\label{Kuznetsovin}
+\frac{1}{\pi}\int_{-\infty}^{\infty}\frac{\sigma_{2ir}(m)\sigma_{2ir}(n)}{(mn)^{ir}|\zeta(1+2ir)|^{2}}\hat{\psi}(r)\, dr 
+\sum_{k=1}^{\infty}a_{k}\sum_{j=1}^{\vartheta(k)}\overline{\rho_{j,k}(m)}\rho_{j,k}(n)\hat{\psi}((1-k)i/2),
\end{equation}
where
\[
\hat{\psi}(r)=\frac{\pi i}{2\sinh (\pi r)}\int_{0}^{\infty}(J_{2ir}(x)-J_{-2ir}(x))\psi(x)x^{-1}\, dx.
\]
\end{lem}
For a proof, see \cite{Motohashi2}, Theorem 2.3. Note that the conditions for $\psi$ can be relaxed. However, the formulation of 
Lemma \ref{Kuzn} suffices for our purposes.

In case $m$ and $n$ are of the opposite sign in the previous lemma, we have the following variation:
\begin{lem}\label{Kuzn2}
Let $\psi(x)$ be a smooth function with compact support in $(0,\infty)$ and $m,n$ be positive integers. Then
\[
\sum_{q=1}^{\infty}q^{-1}S(m,-n;q)\psi\left(\frac{4\pi\sqrt{mn}}{q}\right)=\sum_{j=1}^{\infty}\alpha_{j}\varepsilon_{j}t_{j}(m)t_{j}(n)\hat{\psi}^{-}(\kappa_{j})
\]
\[
+\frac{1}{\pi}\int_{-\infty}^{\infty}\frac{\sigma_{2ir}(m)\sigma_{2ir}(n)}{(mn)^{ir}|\zeta(1+2ir)|^{2}}\hat{\psi}^{-}(r)\, dr,
\]
where
\[
\hat{\psi}^{-}(r)=2\cosh (\pi r)\int_{0}^{\infty}K_{2ir}(x)\psi(x)x^{-1}\, dx,
\]
and $\varepsilon_{j}$ is the parity sign of the $j$th Maass form.
\end{lem}
A proof can be found in \cite{Motohashi2}, Theorem 2.5.

The following lemma gives an estimate for triple sums of Kloosterman sums (see \cite{Iwaniec3}, Theorem 4):
\begin{lem}\label{Iwan}
Let $M,N,Q\geq 1$ and $g(m,n,q)\in C^{2}$ be a weight function
with the properties
\[
\textrm{supp}\; g\subseteq[M,2M]\times[N,2N]\times[Q,2Q]
\]
and, for $0\leq \nu_{1},\nu_{2},\nu_{3}\leq 2$,
\[
\left|\frac{\partial^{\nu_{1}+\nu_{2}+\nu_{3}}}{\partial
m^{\nu_{1}}\partial n^{\nu_{2}}\partial q^{\nu_{3}}}g(m,n,q)\right|\leq
M^{-\nu_{1}}N^{-\nu_{2}}Q^{-\nu_{3}}.
\]
Then for $a_{m},b_{n}\in\mathbb{C}$ we have
\[
\sum_{m\sim M}\sum_{n\sim N}\sum_{q\sim
Q}a_{m}b_{n}g(m,n,q)S(m,\pm n;q)
\]
\[
\ll Q^{1+\varepsilon}(MN)^{1/2}\left(\sum_{m\sim M}
|a_{m}|^{2}\right)^{1/2}\left(\sum_{n\sim N} |b_{n}|^{2}\right)^{1/2}.
\]
\end{lem}

Next we recall a basic inequality in the proof of the
classical large sieve:
\begin{lem}[Sobolev]\label{klassinen}
Let $a\leq u\leq a+\Delta$ for some $a\in \mathbb{R}$, $\Delta\in\mathbb{R_{+}}$, and let the function $f$ be
continuously differentiable on this interval. Then
\[
|f(u)|^{2}\leq \Delta^{-1}\int_{a}^{a+\Delta}|f(x)|^{2}\,dx
+2\left(\int_{a}^{a+\Delta}|f(x)|^{2}\,dx\right)^{1/2}
\]
\begin{equation}\label{Sobol}
\times \left(\int_{a}^{a+\Delta}|f'(x)|^{2}\,dx\right)^{1/2} \ll
\Delta^{-1}\int_{a}^{a+\Delta}\left(|f(x)|^{2}+\Delta^{2}|f'(x)|^{2}
\right)\,dx
\end{equation}
uniformly.
\end{lem}
For a proof, see Montgomery \cite{Montgomery}, Lemma 1.1 applied
to $f^{2}$.

Finally we have by \cite{Kowalski}, pp. 169-170,
\begin{lem}[Duality Principle]\label{dpa}
Let $\Lambda,\Theta\geq 1$. For any complex numbers $b(\lambda)$ and $\phi(\lambda,\theta)$, $1\leq \lambda\leq \Lambda$, $1\leq \theta\leq \Theta$,
\[
\sum_{\theta=1}^{\Theta}\left|\sum_{\lambda=1}^{\Lambda}b(\lambda)\phi(\lambda,\theta)\right|^{2}\leq \sup_{c(\theta)}\sum_{\lambda=1}^{\Lambda}\left|\sum_{\theta=1}^{\Theta}
c(\theta)\phi(\lambda,\theta)\right|^{2} 
\sum_{\lambda=1}^{\Lambda}|b(\lambda)|^{2},
\]
where the supremum over $c(\theta)$ is taken over all complex numbers such that \\ $\sum_{\theta=1}^{\Theta}|c(\theta)|^{2}=1$. 
\end{lem}

Note that the inequality in Lemma \ref{dpa}
holds also when the $\lambda$-sum above is replaced by an integral, if $b(\lambda)$ and $\phi(\lambda,\theta)$ are continuous as functions of $\lambda$. Also
different combinations of several sums and integrals give analogous results.

\section{Proofs of our results}

\subsection{Proof of Lemma \ref{tarkein}}

We shall follow closely the ideas of Jutila's paper \cite{Jutila6}, built on another article by Jutila \cite{Jutila7}. However,
we shall write down the details in order to make this paper as independent as possible and to be reader-friendly when explaining
the needed modifications.

In the sequel, we shall repeatedly use the notation $w_{K}$ for a suitably chosen real-valued smooth weight function, where $0\leq w_{K}(x)\leq 1$ for all $x\in \mathbb{R}$, 
supp $w\subseteq[BK,CK]$ for suitable constants $B,C$, $w_{K}(x)=1$ when $x\asymp K$ and
$w_{K}^{(\nu)}(x)\ll_{\nu}K^{-\nu}$ for each $\nu\geq 0$.

Now by \eqref{tahti}, with the substitution $\eta=x-a/q$,
\[
\sum_{|f|\sim F}\sum_{n\sim N}\left |b_{f}^{\ast}(n)
\right|^{2} =\sum_{|f|\sim F}\sum_{n\sim N}\bigg| \frac{1}{\lambda}\sum_{q=1}^{\infty}w(q)\sum_{\substack{a=1 \\ (a,q)=1}}^{q}\int_{-\infty}^{\infty}
\nu(-\eta)
\]
\[
\times\sum_{m_{1}=1}^{\infty}W_{n}^{0}(m_{1})t(m_{1})
 e\left(m_{1}\left(\frac{a}{q}+\eta\right)\right)\sum_{m_{2}=1}^{\infty}W_{n}(m_{2})t(m_{2})e\left(-m_{2}\left(\frac{a}{q}+\eta\right)\right)
\]
\begin{equation}\label{alku}
\times e\left(-f\left(\frac{a}{q}+\eta\right)\right)\,d\eta\bigg|^{2}.
\end{equation}
Let us first treat the case $f>0$, the case of the opposite sign being analogous.

We shall rewrite the exponential sums $S_{1}=S(W_{n}^{0},a/q+\eta)$ and
$S_{2}=S(W_{n},-a/q-\eta)$ by Lemma \ref{Voronoi}, ending up with integrals of the form
\begin{equation}\label{integraalit}
I_{j}=\int_{0}^{\infty}B_{\pm 2i\kappa}\left(\frac{4\pi\sqrt{m_{j}(z_{j}L+n)}}{q}\right)
W^{j}(z_{j}L+n)e((-1)^{j-1}(z_{j}L+n)\eta)L\, dz_{j},
\end{equation}
where $B_{\nu}$ stands for either a J- or a K-Bessel function, $j=1,2$ and $W^{j}$ equals 
$W_{n}^{0}$, when $j=1$, and $W_{n}$, when $j=2$.

For the J-Bessel function we use the asymptotic expansion
\begin{equation}\label{jiiyksi}
J_{2i\kappa}(x)\sim \left(\frac{2}{\pi x}\right)^{1/2}\cos\left(x-i\kappa\pi-\frac{\pi}{4}\right)
\end{equation}
in case $x\geq N^{\varepsilon}$; see \cite{Lebedev}, Eq. (5.11.6). Note that here and later we use the notation $\sim$ also for this purpose, the meaning being clear from the context. 
Notice that it suffices to study
only the leading term of the expansion, for the others are similar, and the contribution of the rest after sufficiently many terms is negligible.
(See also Convention 2 in \cite{JutilaMotohashi}.) 

In case $x< N^{\varepsilon}$ we apply the formula
\begin{equation}\label{jiikaksi}
J_{\nu}(x)=\frac{(x/2)^{\nu}}{\Gamma(1/2)\Gamma(\nu+1/2)}\int_{0}^{\pi}\cos(x\cos \theta)(\sin \theta)^{2\nu}\, d\theta
\end{equation}
from \cite{Lebedev}, Eq. (5.10.4). Lastly, for the K-Bessel function we always use
\begin{equation}\label{kooyksi}
K_{\nu}(x)=\int_{0}^{\infty}e^{-x\cosh u}\cosh(\nu u)\, du
\end{equation}
from \cite{Lebedev}, Eq. (5.10.23). 

By integrating repeatedly by parts we first notice that 
the contribution of the integral involving the K-Bessel function \eqref{kooyksi} is negligibly small.
In the terms involving the J-Bessel function we divide the new $m_{j}$-sums ($j=1,2$) into two parts $m_{j}< Q^{2}N^{-1+\delta_{1}}$ and $m_{j}\geq Q^{2}N^{-1+\delta_{1}}$, and use 
\eqref{jiikaksi} or \eqref{jiiyksi} from above accordingly. By repeated partial integration we notice that now $m_{j}$ can be
truncated to be $\asymp NQ^{\delta_{1}}$, as essentially was the case in Jutila's paper \cite{Jutila6} also. Note that we have used
the asymmetric choice of our $\nu$-function here.

Next we investigate the integral over $\eta$,
\[
\int_{-\infty}^{\infty}
\nu(-\eta)e(\eta (z_{1}L-z_{2}L-f))\,d\eta,
\]
where $z_{1}$ and $z_{2}$ are the variables of integration from \eqref{integraalit}.
By integrating by parts we notice that we may assume $z_{1}-z_{2}\ll \delta$. 
Hence it is enough to prove that the upper bound of Lemma \ref{tarkein} holds for
 \[
N^{-1}\sup_{\substack{z_{1},z_{2}\asymp 1, \\ z_{1}-z_{2}\ll \delta}}\sum_{f\sim F}\sum_{n\sim N}\bigg| \sum_{q=1}^{\infty}q^{-1}w(q)
\sum_{m_{1}\asymp NQ^{\delta_{1}}}\sum_{m_{2}\asymp NQ^{\delta_{1}}}t(m_{1})t(m_{2})
\]
\[
\times (m_{1}m_{2})^{-1/4}
w_{NQ^{\delta_{1}}}(m_{1})w_{NQ^{\delta_{1}}}(m_{2})
S(f,m_{1}-m_{2};q)
\]
\[
\times \exp\left( \frac{4\pi i}{q}\left(\sqrt{m_{1}(z_{1}L+n)} -\sqrt{m_{2}(z_{2}L+n)}\right)\right)
\bigg|^{2}.
\]

Next we separate the cases $m_{1}>m_{2}$, $m_{1}<m_{2}$ and $m_{1}=m_{2}$. In the last case we write the emerging Ramanujan
sum with its arithmetic formula (\cite{Titchmarsh}, Eq. (1.5.2))
\[
S(f,0;q)=c_{q}(f)=\sum_{d|(q,f)}d\mu\left(\frac{q}{d}\right),
\]
whence
\[
\sum_{q=1}^{\infty}w(q)q^{-1}c_{q}(f)\exp\left(\frac{4\pi i}{q}\left(\sqrt{m(z_{1}L+n)}-\sqrt{m(z_{2}L+n)}\right)\right)
\]
\[
=\sum_{d|f}\sum_{r=1}^{\infty}w(rd)r^{-1}\mu(r)\exp\left(\frac{4\pi i}{rd}\left(\sqrt{m(z_{1}L+n)}-\sqrt{m(z_{2}L+n)}\right)\right).
\]
Hence by trivial estimates the contribution of this case to our sum is
$\ll N^{1+\varepsilon}F$.

We then assume $m_{1}>m_{2}$, and comment the opposite case at the end of the proof.
We make a change of variables $m_{1}=m_{2}+p$, $p>0$, and use the notation $m_{2}=m$ for simplicity. Furthermore we decompose the range of $p$ 
into intervals $p \asymp P$ with $1\leq P\ll NQ^{\delta_{1}}$, the number
of these intervals being $\ll N^{\varepsilon}$, and equip the subsums with suitable weight functions $w_{P}(p)$. As in Jutila's paper \cite{Jutila7} we write
\begin{equation}\label{eeaxaa}
\frac{4\pi i}{q}\left(\sqrt{(m+p)(z_{1}L+n)} -\sqrt{m(z_{2}L+n)}\right)=ixa\left(\frac{m}{p},n\right),
\end{equation}
where
\[
x=\frac{4\pi \sqrt{pf}}{q}\asymp \frac{\sqrt{PF}}{Q}=X
\]
and 
\[
a(y,n)=\frac{z_{1}L+n+yL(z_{1}-z_{2})}{\sqrt{f}\left(\sqrt{(y+1)(z_{1}L+n)}+\sqrt{y(z_{2}L+n)}\right)}.
\]
Now $a(y,n)$ depends on $f$ as well, but since this dependence is less central in the sequel, we shall leave it implicit for the simplicity of the notation.

Following further the ideas of \cite{Jutila6} we next
introduce a new variable $u$ and write the $m$-sum as
\[
\sum_{m=-\infty}^{\infty}G(m)=\frac{1}{U}\sum_{u\asymp U}\sum_{m=-\infty}^{\infty}G(m+u).
\]
In case $P\ll QN^{2\delta_{1}}$, we choose $U\asymp N^{1-2\delta_{1}}$, whence by Taylor's approximation
\[
w_{NQ^{\delta_{1}}}(m+u)w_{NQ^{\delta_{1}}}(m+p+u)\left((m+u)(m+p+u)\right)^{-1/4}e^{ixa((m+u)/p,n)}
\]
\[
\sim w_{NQ^{\delta_{1}}}(m)w_{NQ^{\delta_{1}}}(m+p)\left(m(m+p)\right)^{-1/4}e^{ixa(m/p,n)}.
\]
We add a suitable weight function $w_{F}(f)$ and use Lemma \ref{dpa}, where now our $m$- and $f$-sums correspond to the $\lambda$
and $\theta$-sums, respectively.
We end up to estimate
\[
N^{-1+\varepsilon}U^{-2}Q^{-2}\sup_{c(f)} \sup_{\substack{z_{1},z_{2}\asymp 1, \\ z_{1}-z_{2}\ll \delta}}\sum_{n\sim N}\sum_{m\asymp NQ^{\delta_{1}}}\bigg| \sum_{f\sim F}
c(f)\sum_{q\asymp Q}\sum_{p\asymp P}\sum_{u\asymp U}t(m+u)
\]
\[
\times t(m+p+u)H(f,p,q)S(f,p;q)
\bigg|^{2},
\]
where the supremum over $c(f)$ is taken over all complex numbers depending on $f$ such that $\sum_{f\sim F}|c(f)|^{2}=1$.
Moreover $H$ is a suitable weight function such that for any $0\leq \nu_{1},\nu_{2},\nu_{3}\leq 2$
\[
\left|\frac{\partial^{\nu_{1}+\nu_{2}+\nu_{3}}}{\partial
f^{\nu_{1}}\partial p^{\nu_{2}}\partial q^{\nu_{3}}}H(f,p,q)\right|\leq
F^{-\nu_{1}}P^{-\nu_{2}}Q^{-\nu_{3}}.
\]
Then by Lemma \ref{Iwan} and the earlier estimate \eqref{Tsumma} we end up with the bound $N^{1+\varepsilon}\delta LF$. 

We may therefore assume $P\gg QN^{2\delta_{1}}$. Choosing now
\[
U\asymp\frac{N^{1-\varepsilon}Q}{P}
\]
we notice again by Taylor's approximation that
\[
w_{NQ^{\delta_{1}}}(m+u)w_{NQ^{\delta_{1}}}(m+p+u)\left((m+u)(m+p+u)\right)^{-1/4}e^{ixa((m+u)/p,n)}
\]
\[
\sim w_{NQ^{\delta_{1}}}(m)w_{NQ^{\delta_{1}}}(m+p)\left(m(m+p)\right)^{-1/4}e^{ixa(m/p,n)}.
\]
Furthermore we insert a new weight function $w_{X}\left(\frac{4\pi\sqrt{pf}}{q}\right)$ to our sum, chosen to equal $1$ whenever
$w(q)w_{P}(p)\neq 0$ and $f\sim F$. Now we may write $w(q)=w\left(\frac{4\pi\sqrt{pf}}{x}\right)$ by its Mellin inversion
\[
\frac{1}{2\pi}\int_{-\infty}^{\infty}q^{-c-it}w^{\ast}(c+it)\, dt,
\]
where
\[
w^{\ast}(c+it)=\int_{0}^{\infty}\theta^{c-1+it}w(\theta)\, d\theta,
\]
and we may choose $c=0$. By integrating $w^{\ast}(c+it)$ repeatedly by parts we notice that we may truncate $|t|\ll N^{\varepsilon}$. We move the 
$t$-integral and $w^{\ast}(it)$ out of the square, whence essentially we have replaced $w(q)$ by $w_{X}\left(\frac{4\pi\sqrt{pf}}{q}\right)$. 
This formulation proves to be convenient in the sequel.

Finally we are ready to utilize Kuznetsov's trace formula, Lemma \ref{Kuzn}. We shall first treat the first term on the right hand side of 
\eqref{Kuznetsovin} commenting the other two in the end briefly. We write 
\[
\frac{J_{2ir}(x)-J_{-2ir}(x)}{\sinh(\pi r)}=\frac{4}{\pi i}\int_{0}^{\infty}\cos(x\cosh \xi)\cos(2r\xi)\, d\xi, \;\;\; r\neq 0,x>0,
\]
(see \cite{Lebedev} p.139), and make a change of variables $e^{\xi}=\omega$, whence
\[
\hat{\psi}(r,m/p)=\int_{0}^{\infty}\int_{0}^{\infty}\cos\left(\frac{x}{2}(\omega+\omega^{-1})\right)\cos(2r\ln \omega)w_{X}(x)x^{it}(4\pi\sqrt{pf})^{-it}
\]
\begin{equation}\label{oomega}
\times e^{ixa(m/p,n)}x^{-1}\omega^{-1}\, dx d\omega.
\end{equation}
By repeated integration by parts
we notice that we may assume
$\kappa_{j}\asymp P/Q^{1+\delta_{1}/2}$ and $\omega\asymp \sqrt{P}/(\sqrt{F}Q^{\delta_{1}/2})$ with a negligibly small error.

Next we prepare ourselves to use the duality principle again. As in \cite{Jutila6}, we need to separate
the variable $y=m/p$ from $\hat{\psi}(\kappa_{j},y)$. We therefore write the double sum over $m$ and $p$ as follows:
\[
\sum_{m\asymp NQ^{\delta_{1}}}\sum_{p\asymp P}=\sum_{v=0}^{V}\sum_{\substack{m\asymp NQ^{\delta_{1}}, \; p\asymp P,\\ \frac{m}{p}\in I_{v}}},
\]
where
\[
I_{v}=\left[\frac{BNQ^{\delta_{1}}}{P}+v\Delta,\frac{BNQ^{\delta_{1}}}{P}+(v+1)\Delta\right),
\]
$B>0$ is a suitable constant, $\Delta=N^{1-\varepsilon}Q/P^{2}$ and $V\ll NQ^{\delta_{1}}/(P\Delta)$. The last interval $I_{V}$ may be incomplete. We write
$y_{v}=BNQ^{\delta_{1}}/P+v\Delta$ and express $\hat{\psi}$ again as a Taylor polynomial around $y_{v}$ on each interval $I_{v}$, getting an asymptotic expression
$\hat{\psi}(\kappa_{j},y)\sim \hat{\psi}(\kappa_{j},y_{v})$, when $y\in I_{v}$.
Now we use Lemma \ref{dpa}, where the double sum over $v$ and $\kappa_{j}$ corresponds to the $\lambda$-sum, and the double sum over $f$ and $n$ corresponds
to the $\theta$-sum.
Hence it is enough to prove the bound of Lemma \ref{tarkein} for
 \[
N^{-1}U^{-2}\sup_{c(f,n)}\sup_{\substack{z_{1},z_{2}\asymp 1,\\ z_{1}-z_{2}\ll \delta}}\sup_{|t|\ll N^{\varepsilon}}\sum_{v=0}^{V}\sum_{\kappa_{j}\asymp PQ^{-1-\delta_{1}/2}}
\alpha_{j}\bigg|\sum_{f\sim F}\sum_{n\sim N} c(f,n)t_{j}(f)
\]
\[
\times\hat{\psi}(\kappa_{j},y_{v})\bigg|^{2}
\sum_{v=0}^{V}\sum_{\kappa_{j}\asymp PQ^{-1-\delta_{1}/2}}\alpha_{j}\bigg|\sum_{\substack{m\asymp NQ^{\delta_{1}}, \; p\asymp P,\\ \frac{m}{p}\in I_{v}}}\sum_{u\asymp U}t(m+u)
t(m+u+p)
\]
\[
\times t_{j}(p)w_{P}(p)p^{-it/2}w_{NQ^{\delta_{1}}}(m)w_{NQ^{\delta_{1}}}(m+p)\left(m(m+p)\right)^{-1/4}\bigg|^{2},
\]
where again the supremum over $c(f,n)$ is taken over all complex numbers such that $\sum_{f\sim F}\sum_{n\sim N}|c(f,n)|^{2}=1$.

We next apply the spectral large sieve to the latter $\kappa_{j}$-sum, take the $m$-sum out of the square by Cauchy's inequality and use
again \eqref{Tsumma}.
With the first $\kappa_{j}$-sum we face the problem of $\hat{\psi}(\kappa_{j},y_{v})$ depending
on $\kappa_{j}$, which we overcome by use of Sobolev's Lemma \ref{klassinen} using Jutila's paper \cite{Jutila1}, p. 454,  as a model:   
The range $\kappa_{j}\asymp P/Q^{1+\delta_{1}/2}$ is
split up into segments of length $1$, whence $\hat{\psi}(\kappa_{j},y_{v})$ remains essentially stationary as $\kappa_{j}$ runs
over a segment. In this way, the second term on the right hand side of
\eqref{Sobol} will be comparable to the first. 
Hence we divide the $\kappa_{j}$-sum into
subsums of length $1$, and apply Lemma \ref{klassinen} to each subsum. Next we apply the spectral large sieve to each subsum over
$\kappa_{j}$, and finally add the results
together. This leads us to the bound
\[
N^{\varepsilon}PQ^{-1}\left(\left(\frac{P}{Q}\right)^{2}+P\right)\left(\frac{P}{Q}+F\right)\sup_{0\leq v\leq V}\sup_{\substack{c(f,n), t,\\ z_{1},z_{2}}}\int\sum_{f}\bigg(\bigg| \sum_{n}c(f,n)
\]
\[
\times\hat{\psi}(r,y_{v})\bigg|^{2}
+\bigg|\sum_{n}c(f,n)
\frac{\partial}{\partial r}\hat{\psi}(r,y_{v})\bigg|^{2}\bigg)
w_{P/Q^{1+\delta_{1}/2}}(r)\, dr.
\]

We proceed to utilize the averaging over $n$. 
By repeated integration by parts over $x$ in \eqref{oomega} we notice that we may assume $\omega-2a(y_{v},n)\ll N^{2\delta_{1}}X^{-1}$, otherwise the $x$-integral is negligibly small.
Lastly, we open the squares, whence 
the $\omega$-integral from \eqref{oomega} produces two integrals, say, over variables $\omega_{1}$ and $\omega_{2}$. By integrating repeatedly by parts with respect to
$r$ we notice that we may further truncate $\omega_{1}-\omega_{2}\ll N^{2\delta_{1}}X^{-1}$. Therefore also $a(y_{v},n_{1})-a(y_{v},n_{2})\ll N^{2\delta_{1}}X^{-1}$, whence
by the mean value theorem $n_{1}-n_{2}\ll N^{1+3\delta_{1}}Q/P$. 

Hence we finally end up with the upper bound
\[
N^{\varepsilon}\left(\left(\frac{P}{Q}\right)^{2}+P\right)\left(\frac{P}{Q}+F\right)\sup_{c(f,n)}
\sum_{f\sim F}\sum_{n\sim N}|c(f,n)|^{2}\sum_{|n-n'|\ll N^{1+\varepsilon}QP^{-1}}1,
\]
and conclude with the desired result.

Now  the treatment of the second term in Lemma \ref{Kuzn} is completely analogous, as we have Lemma
\ref{jatkuvasss} to mimic the spectral large sieve.

The third term produces only a negligible contribution, as can be seen directly for example by the equation
\[
\frac{\pi i}{2\sinh\left((1-k)\frac{\pi i}{2}\right)}\left(J_{k-1}(x)-J_{-(k-1)}(x)\right)
\]
\[
=(-1)^{1+k/2}\int_{0}^{\pi}\sin(x\sin \theta)\sin((1-k)\theta)\, d\theta, \;\;\; x>0, k\equiv 0\; (\textrm{mod}\; 2),
\]
(see \cite{Lebedev}, Eq. (5.10.8)), and by the repeated partial integration, using \eqref{aakoo}.

In case $m_{1}<m_{2}$ we use the notation $m_{1}=m_{2}+p$, whence again $p>0$. The deduction is analogous to that above,
except that instead of Lemma \ref{Kuzn} we use Lemma \ref{Kuzn2} with
\[
\cosh(\pi r)K_{2ir}(x)=\int_{0}^{\infty}\cos(x\sinh \xi)\cos(2r\xi)\, d\xi, \;\; x>0,
\]
from \cite{Watson}, Eq. (13), p.183.

\subsection{Proof of Lemma \ref{nolla}}

By \eqref{beexii} and Cauchy's inequality
\[
\sum_{n\sim N}|b_{0}^{\ast}(n)|^{2}\ll \sum_{n}\int_{0}^{1}|\chi^{\ast}(x)|^{2}|S(W_{n}^{0},x)|^{2}\, dx
\int_{0}^{1}|S(W_{n},-x)|^{2}\, dx.
\]
Now we use the fact that $\chi^{\ast}(x)\ll 1$, open the rest of the squares, integrate and end our proof with the help of \eqref{Hoff}.

\subsection{Proof of Lemma \ref{johtopaatos}}

We shall follow precisely the steps of the proof of Theorem 2 in \cite{Jutila6}.
Hence, by \cite{Jutila6}, Eq. (1.7), (as a direct consequence of \eqref{khii})
\[
\sum_{n\sim N}\left |b_{f}(n)\right|^{2} \ll \sum_{n\sim N}\left |b_{f}^{\ast}(n)
\right|^{2} +\sum_{n\sim N}\left |\frac{1}{\lambda}\sum_{d=1}^{\infty}d\sum_{m\neq 0}a_{dm}b_{f+dm}(n)\sum_{r=1}^{\infty}w(dr)\mu(r)
\right|^{2}.
\]
By \eqref{aFourier} the summation over $d$ and $m$ such that $|dm|>\delta L$ yields only negligibly small contribution.
Therefore
\[
\sum_{n\sim N}\left |\frac{1}{\lambda}\sum_{d=1}^{\infty}d\sum_{m\neq 0}a_{dm}b_{f+dm}(n)\sum_{r=1}^{\infty}w(dr)\mu(r)
\right|^{2}\ll Q^{-2}N^{\varepsilon}\sum_{n\sim N}\left |\sum_{1\leq |\xi|\leq \delta L}|b_{f+\xi}(n)|
\right|^{2}
\]
\[
\ll (\delta L)^{-1}N^{\varepsilon}\sum_{n\sim N}\sum_{-\delta L\ll f\ll \delta L}|b_{f}(n)|^{2}.
\]

Now by Lemma \ref{uusiarvio} and Cauchy's inequality
\[
\sum_{-\delta L\ll f\ll \delta L}\sum_{n\sim N}\left |b_{f}(n)
\right|^{2} \ll N^{-A}+Q^{\varepsilon}\max_{|\xi|\ll Q}\sum_{-\delta L\ll f\ll \delta L}\sum_{n\sim N}|b_{f+\xi}^{\ast}(n)|^{2}
\]
\[
\ll N^{-A}+Q^{\varepsilon}\sum_{-\delta L\ll f\ll \delta L}\sum_{n\sim N}|b_{f}^{\ast}(n)|^{2},
\]
which is by Lemmas \ref{tarkein} and \ref{nolla} 
\[
\ll N^{3+\varepsilon}(\delta L)^{-2}+N^{2+\varepsilon}+N^{1+\varepsilon}(\delta L)^{2}.
\]

Therefore finally
\[
 \sum_{n\sim N}\left |b_{f}(n)
\right|^{2}\ll \sum_{n\sim N}\left |b_{f}^{\ast}(n)
\right|^{2}+ N^{3+\varepsilon}(\delta L)^{-3}+N^{2+\varepsilon}(\delta L)^{-1}+N^{1+\varepsilon}\delta L.
\]

\subsection{Proof of Theorem \ref{thm1}}

Let us insert another weight function
$W_{n}^{0}(n+l+f)$ to the $l$-sum, where
we let $0\leq W_{n}^{0}(x)\leq 1$ stand for a real-valued smooth weight function supported on the interval 
$[(B-\delta)L+n,(C+(D+1)\delta)L+n]$, with
$D$ a suitable large constant, $W_{n}^{0}(x)=1$ on
$[BL+n,(C+D\delta)L+n]$ and
$(W_{n}^{0})^{(\nu)}(x)\ll_{\nu}(\delta L)^{-\nu}$ for each $\nu\geq 0$ and $x\in \mathbb{R}$. 
Note that now $W_{n}^{0}(n+l+f)=1$ whenever $W_{n}(n+l)\neq 0$. 

Then the theorem is a direct consequence of Lemmas \ref{tarkein} and \ref{johtopaatos}. 
Note that we may add the interval $1\leq L\ll N^{\varepsilon}$ by trivial estimations.

\subsection{Proof of Theorem \ref{thm2}}

First we notice that if $N^{2}\gg FL^{4}$, then the trivial estimates are enough. Hence we may assume that $N^{2}\ll FL^{4}$.

Let $L^{-1+\varepsilon}\leq U_{1}\leq 1/4$ be a quantity to be fixed later, and suppose that $F\ll U_{1}L$. We insert a set of $\asymp\log(1/U_{1})$ real-valued smooth weight functions
$g_{\delta}\left(\frac{l}{L}\right)$ to the $l$-sum so that
their sum produces an approximation of the characteristic function of the interval $[1,2]$ with an error of size $\ll U_{1}$
and their supports 
widen step by step by factors 2 when we move away from the end points 1 and 
2. To be precise, let the first weight function $g_{U_{1}}(x)$ to be supported on the interval $[1,1+4U_{1}]$, and $g_{U_{1}}(x)=1$ when $x\in[1+U_{1},1+2U_{1}]$. 
Then $g_{2U_{1}}(x)$ is supported on the interval $[1+2U_{1},1+10U_{1}]$, and $g_{2U_{1}}(x)=1$ when $x\in[1+4U_{1},1+6U_{1}]$ and so on, and the slopes of the weight functions cancel out each other. 
Hence for $U_{1}\leq \delta\leq 1/4$ the function $g_{\delta}(x)$ is a compactly supported function on some interval $[B,C]$ of length $\asymp\delta$ contained in $[1,2]$. Moreover,  $g_{\delta}(x)=1$ on
an interval of length $\asymp\delta$ and
$g_{\delta}^{(\nu)}(x)\ll_{\nu}\delta^{-\nu}$ for each $\nu\geq 0$ and $x\in \mathbb{R}$. 

Therefore by \eqref{Hoff} 
\[
\sum_{f\sim F}\sum_{n\sim N}\left|\sum_{l\sim L}t(n+l)t(n+l+f)
\right|^{2}
\]
\[
\ll L^{\varepsilon}\sup_{U_{1}\leq \delta\leq\frac{1}{4}}\sum_{f\sim F}\sum_{n\sim N}\bigg |\sum_{l\sim L}t(n+l)t(n+l+f)g_{\delta}\left(\frac{l}{L}\right)\bigg|^{2}
+N^{1+\varepsilon}FU_{1}^{2}L^{2}.
\]
Denoting $g_{\delta}(x)$ by $W_{n}(xL+n)$ we conclude the desired result by
Theorem \ref{thm1}, choosing 
\[
U_{1}\asymp\frac{N^{1/2}}{F^{1/4}L}.
\]
Note that  $L^{-1+\varepsilon}\leq U_{1}\leq 1/4$ is satisfied by our assumption $N^{2}\ll FL^{4}$. Furthermore, the condition $F\ll U_{1}L$ yields the requirement
$F\ll N^{2/5}$.

\subsection{Proof of Theorem \ref{thm4}}

The deduction is analogous to that above. The Vorono\u\i \, type formula can be found in \cite{Jutila0}, Theorem 1.7, and instead of \eqref{Tsumma} we now use
\eqref{aasumma}.

\subsection*{Acknowledgment}
The author wishes to thank Professor Jutila for all the support and helpful advices. Especially
he pointed out the result presented in Lemma \ref{uusiarvio} from his preprint \cite{Jutila10} to overcome problems in estimating the 
difference $b_{f}-b_{f}^{\ast}$, and gave his permission to publish this result in this paper. Also the valuable comments by the anonymous reviewer for clarifying and simplifying
the arguments are gratefully acknowledged. The help from these two truly improved the quality of this paper.

\vspace{1cm}

{\sc Eeva Suvitie}

{\sc Department of Mathematics}

{\sc FI-20014 University of Turku} 

{\sc Finland}

{\sc e-mail: eeva.suvitie@utu.fi}

\end{document}